\documentclass{amsart}

\usepackage{pstricks, pst-tree, pst-node} 
\usepackage[latin1]{inputenc}
\usepackage{amssymb}
\usepackage{xfrac}     
\usepackage{lmodern}
\usepackage{stmaryrd}
\usepackage{graphicx}
\usepackage{amsmath}  
\usepackage{amsthm}
\usepackage{wrapfig} 
\usepackage{color}  
\usepackage{microtype} 
\usepackage{setspace}
\usepackage{amsfonts}
\usepackage{wasysym}
\usepackage{cancel}
\usepackage{verbatim}
\usepackage{dsfont}
\usepackage{paralist}
\usepackage{mathtools}

\DeclareMathOperator{\spann}{span}

\DeclareMathOperator{\p}{\varphi}

\DeclareMathOperator{\FF}{\mathbb{F}}
\DeclareMathOperator{\lr}{\langle}
\DeclareMathOperator{\rr}{\rangle}

\DeclareMathOperator{\ZZ}{\mathbb{Z}}

\DeclareMathOperator{\NN}{\mathbb{N}}

\DeclareMathOperator{\lrp}{\langle\!\langle}
\DeclareMathOperator{\rrr}{\rangle\!\rangle}

\DeclareMathOperator{\bb}{\mathfrak{b}}
\DeclareMathOperator{\pb}{\varphi_{\mathfrak{b}}}

\DeclareMathOperator{\w}{\wedge}
\DeclareMathOperator{\de}{d\!}

\DeclareMathOperator{\pB}{\mathcal B}
\DeclareMathOperator{\coker}{coker}

\DeclareMathOperator{\ndeg}{ndeg}

\DeclareMathOperator{\pdeg}{\textit{p}\hspace{1.5pt}-deg}

\DeclareMathOperator{\ann}{Ann}

\newcommand{\be}{\hfill$\square$}

\newtheoremstyle{plain2}
  {10pt}   
  {10pt}   
  {\itshape}  
  {0pt}       
  {\bfseries} 
  {}         
  {5pt plus 1pt minus 1pt} 
  {}          
  
 \newtheoremstyle{beweis}
  {10pt}   
  {10pt}   
  {\normalfont}  
  {0pt}       
  {\bfseries} 
  {:}         
  {5pt plus 1pt minus 1pt} 
  {}          

\newtheoremstyle{definition2}
  {10pt}   
  {10pt}   
  {\normalfont}  
  {0pt}       
  {\bfseries} 
  {}         
  {5pt plus 1pt minus 1pt} 
  {}          

\theoremstyle{plain2}
\newtheorem{satz}{Satz}[section]
\newtheorem{lem}[satz]{Lemma}
\newtheorem{pro}[satz]{Proposition}
\newtheorem{theo}[satz]{Theorem}

\theoremstyle{definition2}
\newtheorem{defi}[satz]{Definition}
\newtheorem{rem}[satz]{Remark}

\theoremstyle{beweis}
\newtheorem*{prf}{Proof}

\begin{document}

\title[Annihilators and function field extensions of p-forms]{Function field extensions and annihilators of differential forms in characteristic $p$ and bilinear forms in characteristic $2$}

\author{Marco Sobiech}
\address{Fakult\"at f\"ur Mathematik, Technische Universit\"at Dortmund, D-44221 Dortmund, Germany}
\email{marco.sobiech@tu-dortmund.de}
\date{\today}

\begin{abstract}
Let $F$ be a field of characteristic $p>0$ and let $\Omega^n(F)$ be the $F$-vector space of $n$-differential forms over $F$. In this work we will study the 
behaviour of $\Omega^n(F)$ under iterated function field extensions of $p$-forms. We will use results from a previous work to rewrite the kernel 
of the restriction map $\Omega^n(F) \to \Omega^n(F(\p_1,\ldots,\p_r))$ 
to a set of differential forms annihilated by specific forms given by the norm fields of  $\p_1,\ldots,\p_r$.
\end{abstract}

\subjclass[2010]{Primary 
12F99, 12H05, 13N05, 14F99 }

\keywords{differential forms; Artin-Schreier-map; annihilator}

\maketitle

\begin{section}{Introduction}

Let $F$ be a field of characteristic two and 
let $W(F)$ resp. $W_q(F)$ denote the 
bilinear Witt ring resp. the quadratic Witt group of $F$. These two 
algebraic structures heavily depend on the base field $F$, so 
it is a natural question to ask how these structures behave under field 
extensions. For separable extensions $S/F$, the bilinear Witt kernel 
$W(S/F)$ is trivial, since any $2$-basis of $F$ is still 
a $2$-basis of $S$. The bilinear Witt kernel of an arbitrary purely inseparable 
extension $E/F$ however is noticeably more difficult to compute. If 
$E/F$ is a modular purely inseparable extension, i.e. 
$E=F\left( \sqrt[2^{m_1}]{a_1},\ldots, \sqrt[2^{m_r}]{a_r}\right)$ for 
some $2$-independent $\{a_1,\ldots,a_r\}\subset F$, we have full control 
over the $2$-basis of $E$ by starting with a suitable $2$-basis over $F$  
and hence can give a complete description of the generators of the 
kernel $W(E/F)$. If 
we drop the modularity of the extension, things become more difficult rather 
quickly. In a previous work, we introduced the theory of 
annihilators of differential forms so that we were able to give 
a generating system of the bilinear Witt kernel
 for an arbitrary two fold purely inseparable extension. 

In this paper, we will now study the bilinear Witt kernel 
for non-algebraic field extensions. The most prominent case of such an 
extension is probably the case of a function field extension of a 
quadratic form. Let $\pi=\lrp a_1,\ldots,a_r\rrr$ be the 
polar form of an anisotropic bilinear Pfister 
form over $F$ (these forms are usually called quasi-Pfister forms). 
Then both kernels $W(F(\pi)/F)$ and $W_q(F(\pi)/F)$ are 
known and given by 
\begin{align*}
W(F(\pi)/F) &= \left[ \lrp y_1,\ldots,y_r\rrr \mid  y_1,\ldots,y_r\in F^2(a_1,\ldots,a_r)^* \right] 
\end{align*}
and
\begin{align*}
W_q(F(\pi)/F)&= \pi \otimes W_q(F).
\end{align*}
Note that $F^2(a_1,\ldots,a_r)$ is precisely the norm field of $\pi$. The case of not having the function field of one but of $r\in\NN$ 
quasi-Pfister forms $\pi_1, \ldots,\pi_r$ is studied in this paper. 
For this we will first generalize our problem in the following way: 
Let $F$ now be a field of characteristic $p$. We will 
not study function field extensions of bilinear forms but of 
quasilinear $p$-forms over $F$. We will further translate our problem 
to the theory of differential forms over $F$ so that we are able to 
use the theory of annihilators of differential forms introduced 
in \cite{q58}. With this we will rewrite our studied kernels as 
specific annihilators using the norm fields of  quasilinear $p$-forms 
$\p_1, \ldots,\p_r$. Then we can give generating systems for the kernels 
in every case the annihilator is known. We then come back to the theory of 
bilinear forms of fields of characteristic two via the famous Theorem 
of Kato, which not only answers the Milnor 
conjecture in characteristic two, 
but also connects bilinear forms (and quadratic forms) to differential forms 
in the case $p=2$. 

In \cite{q26} Aravire, Laghribi and O'Ryan studied 
the quadratic Witt kernel of the extension $L(\pi)/F$ with 
$L=F(\sqrt{b_1},\ldots, \sqrt{b_t})$, $\pi=\lrp a_1,\ldots,a_r\rrr$ being 
an anisotropic quasi-Pfister form over $F$ and the additional assumption that 
$\pi$ remains anisotropic over $L$ (which is equivalent to 
$\{a_1,\ldots,a_r,b_1,\ldots,b_t\}$ being $2$-independent over $F$). 
We will close 
our work by giving a more general bilinear counterpart of this kernel, 
since we are able change $\pi$ to an arbitrary bilinear form over $F$ 
and can drop the assumption of $2$-independence.

\end{section}

\begin{section}{A short introduction to $p$-forms}\label{pform}

We will start this paper by briefly recalling the basic facts of quasilinear $p$-forms, 
or $p$-forms for short, over a field $F$ of characteristic $p>0$. These forms were introduced by Hoffmann in 
2004 in \cite{q11} as a generalization of totally singular quadratic forms over fields of characteristic two. We refer to \cite{q11} or \cite{q46} 
for a more detailed introduction to $p$-forms.

Now let $F$ always be a field of characteristic $p>0$, if not stated otherwise and let $V$ be a vector space over $F$. A map $\p: V \to F$ with the properties
\begin{itemize}
	\item $\p(\lambda v)=\lambda^p \p(v)$ for all $\lambda \in F$, $v \in V$
	\item $\p(u+v)=\p(u)+\p(v)$ for all $u,v\in V$
\end{itemize}
is called a $p$-form on $V$ (we will drop the $V$ if it is not necessary). We set $\dim \p= \dim V$ and we 
will focus on finite dimensional $p$-forms for the remainder of the paper. For a field extension $E/F$, by $\p_F$ we will denote the $p$-Form $\p$ 
viewed as a $p$-form over the vector space $V\otimes_F E$.
We further define the set of elements represented by $F$ as 
$D_F(\p):=\{ \p(v) \neq 0 \mid v \in V\}$ and $D_F^0(\p):=D_F(\p)\cup \{0\}$. Let $\p$ resp. $\psi$ be a $p$-form on $V$ resp. $U$. The forms are called isometric, $\p \cong \psi$ for short, if there exists a linear transformation $L: V \to U$ with $\p(v)=\psi(L(v))$ for all $v \in V$. As usually we call $\p$ isotropic, if there is $0\neq v \in V$ 
with $\p(v)=0$ and anisotropic otherwise. For a given basis $v_1,\ldots,v_n$ of $V$ with $a_i=\p(v_i)$ for $i=1,\ldots, n$, we denote the $p$-form 
$\p$ by $\lr a_1,\ldots,a_n\rr_p$. Note that now the $p$-form can be identified with the polynomial $\sum_{i=1}^n a_iX_i^p$ and every polynomial of this 
type defines a $p$-form after fixing a basis. 

For two $p$-forms $\p=\lr a_1,\ldots,a_n\rr_p$, $\psi=\lr b_1,\ldots, b_m\rr_p$ ($a_i,b_j \in F$) and $k\in \NN$ we set
\begin{itemize}
	\item $\p \perp \psi = \lr a_1,\ldots,a_n, b_1,\ldots, b_m\rr_p$,
	\item $k \times \p = \perp_{i=1}^k \p$,
	\item $\p \otimes \psi = \perp_{i,j} \lr a_ib_j\rr_p$.
\end{itemize}
These operations correspond with the usual orthogonal sum and tensor product of vector spaces and are well defined. Note that there is no 
law of cancellation for $p$-forms since we have $\lr 1,1\rr_p\cong \lr 1,0\rr_p$ over $\FF_p$, but $\lr 1\rr_p \not\cong \lr 0\rr_p$. However the 
 question of isometry of $p$-forms may be reduced to basic linear algebra. Since $F^p$ is a subfield of $F$,  
$D_F^0(\lr a_1,\ldots,a_n\rr_p)=\spann_{F^p}(a_1,\ldots,a_n)$ is a subspace of $V$ and we get the following 

\begin{pro}{\cite[Proposition 2.6]{q11}}\label{21}
Let $\p=\lr a_1,\ldots,a_n\rr_p$ and $\psi=\lr b_1,\ldots,b_m\rr_p$ be $p$-forms.
\begin{enumerate}[(a)]
	\item If $\dim_{F^p}D_F^0(\p)=k<n$, choose an $F^p$-basis $c_1,\ldots,c_k$ of $D_F^0(\p)$. Then 
	$\p \cong \lr c_1,\ldots,c_k\rr_p \perp (n-k)\times \lr 0\rr_p$. 
	
	In particular we have $\p \cong \p_{an} \perp (\dim \p - \dim \p_{an})\times \lr 0 \rr_p$, with $\p_{an}$ being the anisotropic part 
	of $\p$, which is unique up to isometry.
	
	\item We have $\p \cong \psi$ if and only if $n=m$ and $D_F^0(\p)=D_F^0(\psi)$.
\end{enumerate}
\end{pro}
 The number $\dim \p - \dim \p_{an}$ in Proposition \ref{21} is called the defect of $\p$ and is noted as $i_d(\p)$.

Again let $\p$ be a $p$-form. We define 
$$N_F(\p):=F^p\left( \frac ab \mid a,b \in D_F(\p) \right)$$
as the norm field of $\p$ and we set $\ndeg_F(\p)=[N_F(\p):F^p]$ as the norm degree of $\p$ over $F$. Note that the 
norm degree of a $p$-form is always a power of $p$ since $N_F(\p)/F^p$ is a purely inseparable extension 
of exponent one. Also note that $\p$ and $x\p$ for $x \in F^*$ share the same norm field by definition. It is commonly known that 
for a $p$-form $\p= \lr a_0,a_1,\ldots,a_n\rr_p$ with $n\in\NN$ and $a_0\neq 0$, the norm field may be 
written as $N_F(\p)=F^p\left(\frac{a_1}{a_0},\ldots, \frac{a_n}{a_0}\right)$.

We will now introduce the so called function field extension of $p$-forms, which we will use frequently in the upcoming chapters. To begin, we need the 
following

\begin{lem}{\cite[Lemma 7.1]{q11}}\label{22}
Let $\p=\lr a_1,\ldots,a_n\rr_p$ be a non-zero $p$-form. Then the corresponding 
polynomial $\p(X_1,\ldots,X_n)$ is irreducible in $F[X_1,\ldots,X_n]$ iff $\ndeg_F(\p)>1$.
\end{lem}

Now again let $\p=\lr a_1,\ldots,a_n\rr_p$ be a non-zero $p$-form and set $X:=(X_1,\ldots,X_n)$. The function field of $\p$, denoted by $F(\p)$, 
is defined as 
$$F(\p):= \begin{cases} \text{Quot}(F[X]/(\p(X)), \ \text{if } \ndeg_F(\p)>1, \\ F(X), \ \text{if } \ndeg_F(\p)=1. \end{cases}$$
To get a better understanding of these kind of fields, we will close this section with some properties of these extensions.

\begin{lem}{\cite[Remark 7.4 + 7.6]{q11}}\label{23}
Let $\p$ and $\psi$ be non-zero $p$-forms.
\begin{enumerate}[(a)]
	\item For $\dim \p \geq 2$, the form $\p_{F(\p)}$ is isotropic.
	\item For $x\in F^*$ we have $F(\p)=F(x\p)$.
	\item For $\p=\psi \perp t\times \lr 0\rr_p$, we have $F(\p)=F(\psi)(Y_1,\ldots,Y_t)$ for some variables $Y_1,\ldots,Y_t$. 
	In particular $F(\p)/F(\p_{an})$ is a purely transcendental extension of transcendence degree $i_d(\p)$.
	\item If $\p=\lr 1,a_1,\ldots,a_n\rr_p$ for $a_1,\ldots,a_n \in F$ with $\ndeg_F(\p)>1$, then 
	$$ F(\p)=F(X_1,\ldots,X_n)\left(\sqrt[p]{\sum_{i=1}^n a_i X^p}\, \right)=F(X_1,\ldots,X_n)\left(\sum_{i=1}^n\sqrt[p]{a_i}X \right).$$
\end{enumerate}
\end{lem}

So by part (d) of Lemma \ref{23}, any non purely transcendental function field extension of a $p$-form can be obtained by a 
purely transcendental extension of transcendence degree $\dim \p-1$ followed by a purely inseparable extension of degree $p$.

\end{section}

\begin{section}{A short introduction to the theory of differential forms}

\begin{subsection}{Basic definitions}

We refer to \cite{q22} or \cite{q1} (for the case $p=2$) for any unmentioned terminology and for a deeper introduction to this topic. 

For the study of differential forms, one of our main tools is the so called $p$-independence. Recall that $F^p$ is a subfield of $F$ and a 
subset $\mathcal A \subset F$ is called $p$-independent, if for every finite subset $\{a_1,\ldots,a_k\} \subset \mathcal A$ we have 
$[F^p(a_1,\ldots,a_k):F^p]=p^k$. The set $\mathcal A$ is called $p$-basis of $F$, if $\mathcal A$ is $p$-independent with $F^p(\mathcal A)=F$. For 
any $\mathcal C \subset F$ we call $\pdeg_F(\mathcal C):=\log_p[F^p(\mathcal C):F^p]$ the $p$-degree of $\mathcal C$ over $F$, if $[F^p(\mathcal C):F^p]$ is 
finite. Otherwise we say that $\mathcal C$ has infinite $p$-degree. Now let us recall some facts about differential forms.

The set $\Omega^1(F)$ of 1-differential forms is defined as the $F$-span of the symbols $\de a$, $a\in F$ with respect to the properties 
$\de\left(ab\right)=a\de b + b \de a $ and $\de\left(a+b\right)=\de a + \de b$ for all $a,b\in F$. Note that $\de: F \to \Omega^1(F), x \mapsto \de x$ is an $F^p$-derivation. 
Using the $n$-th exterior power with $n\in \NN$, we obtain the space of $n$-differential forms by $\bigwedge^n \Omega^1(F)=:\Omega^n(F)$. Thus 
$\Omega^n(F)$ is an $F$-vector space generated by the elementary wedge products $\de a_1 \w \ldots \w \de a_n$ with $a_1,\ldots,a_n \in F$. A scaled generator 
of the type $\frac{\de a_1}{a_1} \w \ldots \w\frac{\de a_n}{a_n}$ with $a_1,\ldots,a_n \in F^*$ is called a logarithmic differential form.  For completeness, we 
set $\Omega^0(F)=F$ and $\Omega^z(F)=\{0\}$ for $z<0$.

By abuse of notation, the map $\text{d}$ can be extended to an $F^p$-linear map $\de :\Omega^{n-1}(F) \to \Omega^{n}(F)$ defined on additive generators by 
$\de \left(x  \de a_1 \w \ldots \w \de a_{n-1}\right)= \de x \w \de a_1 \w \ldots \w \de a_{n-1}$. 
The image of $\text{d}$, denoted by $\de \Omega^{n-1}(F)$, is called the set of exact forms in $\Omega^{n}(F)$. Note that $\de \circ \de =0$. 

The connection between differential forms and $p$-independence can now be seen in the following lemma.

\begin{lem}\label{31}
	\begin{enumerate}[(a)]
		\item Let $a_1,\ldots,a_n\in F$. Then the following statements are equivalent
			\begin{enumerate}[(i)]
				\item $\{a_1,\ldots,a_n\}$ is $p$-independent.
				\item $\de a_1,\ldots,\de a_n$ are $F$-linearly independent.
				\item $\de a_1 \w \ldots \w \de a_n \neq 0$ in $\Omega^n(F)$.
			\end{enumerate}
		\item Let (a) be true for $a_1,\ldots,a_n\in F$ and let $b_1,\ldots,b_n\in F$. Then the following statements are equivalent
		\begin{enumerate}[(i)]
			\item $F^p(a_1,\ldots,a_n)=F^p(b_1,\ldots,b_n)$.
			\item $\spann_F(\de a_1,\ldots, \de a_n)=\spann_F(\de b_1,\ldots,\de b_n)$.
			\item $\de a_1\w\ldots\w\de a_n= x \de b_1 \w \ldots \w \de b_n$ for some $x \in F^*$.
		\end{enumerate}
		
		\item $\de a \in \spann_F(\de b_1,\ldots,\de b_n)$ if and only if  
$a \in F^p(b_1,\ldots,b_n)$.
	\end{enumerate}
	\end{lem}

Now let $\mathcal B=\{b_i \mid i \in I\}$ be a $p$-basis of $F$ for some well-ordered index set $I$. By Lemma \ref{31} we see 
that $\Omega^n(F)=\{0\}$ iff $n> |\mathcal B|$ and by using basic linear algebra, we see that every $p$-independent set can be 
extended to a $p$-basis of $F$. Furthermore the set $\{ \de b_i \mid i \in I\}$ is an $F$-basis of $\Omega^1(F)$, thus 
\begin{equation*}
\bigwedge\nolimits^n_{\pB}:=\{ \de b_{i_1} \w \ldots \w \de b_{i_n} \mid i_1,\ldots,i_n\in I, i_1 < \ldots < i_n\}
\end{equation*} 
is an $F$-basis of $\Omega^n(F)$.

Now let $\wp$ denote the usual Artin-Schreier map $\wp: F \to F, x \mapsto x^p-x$. Then $\wp$ can be extended to the $\Omega^n(F)$ defined on 
logarithmic generators by 
\begin{align*}
\wp : \Omega^n(F) &\to \Omega^n(F) / \de \Omega^{n-1}(F)\ , \\
x \frac {\de a_1}{a_1}\w \ldots \w \frac{\de a_n}{a_n} &\mapsto
(x^p-x) \frac {\de a_1}{a_1}\w \ldots \w \frac{\de a_n}{a_n} \mod \de \Omega^{n-1}(F).
\end{align*}
Due to the Theorem of Kato \cite[Theorem on p.494]{q9} (see also section \ref{bf}), we are in particular interested in the 
kernel and the cokernel of $\wp$, which we will denote by $\nu_n(F):= \ker \wp$ and by $H^{n+1}_p(F):= \coker \wp$. Kato also 
showed in \cite{q4} that $\nu_n(F)$ is additively generated by the logarithmic differential forms 
$\frac{\de a_1}{a_1} \w \ldots \w\frac{\de a_n}{a_n}$ with $a_1,\ldots,a_n \in F^*$.

Let $K/F$ be a field extension. A differential form $\omega\in \Omega^n(F)$ viewed as an element of $\Omega^n(K)$ will be 
shortly denoted by $\omega_K$. The main goal of this paper is to describe which forms of $\Omega^n(F)$ become the zero form over $K$,
so we are interested in the kernels of the restriction map induced by $F \to K$ which we will denote by 
\begin{align*}
\Omega^n(K/F)&:= \{ \omega \in \Omega^n(F) \mid \omega_K =0 \text{ in } \Omega^n(K) \}, \\
\nu_n(K/F)&:= \{ \chi \in \nu_n(F) \mid \chi_K =0 \text{ in } \nu_n(K) \}= \Omega^m(K/F) \cap \nu_n(F).
\end{align*}
We will call $\Omega^n(K/F)$ the $\Omega$-kernel of $K/F$ and $\nu_n(K/F)$ 
the $\nu$-kernel of $K/F$.

In many of our upcoming arguments, the kernel of a purely transcendental extension plays an important role. For easy reference, we will 
state it here as

\begin{pro}\cite[Lemma 2.2]{q1}\label{315}
Let $K/F$ be a purely transcendental extension. Then
$$ \Omega^n(K/F) = \{0\}\ \text{ and }\ \nu_n(K/F)=\{0\}.$$
\end{pro}

\end{subsection}

\begin{subsection}{A recap of annihilators of differential forms}

To get a hand on the kernels $\Omega^n(K/F)$ and $\nu_n(K/F)$ for some non algebraic field extensions $K/F$ we will use annihilators of 
differential forms. These were studied in an earlier paper \cite{q58} and we refer to this source for further details and specific proofs. 

\begin{defi}\label{32}
For $n,r \in \NN$ and a non-empty set $U \subset \Omega^r(F)$, we set 
\begin{align*}
\ann \Omega^n_F(U):&= \{ \omega\in \Omega^n(F)  \mid \omega\w u =0 \in \Omega^{n+r}(F) \text{ for all } u\in U  \} \\
\ann \nu^n_F(U):&= \{ \chi\in \nu_n(F)  \mid \chi\w u =0\in\Omega^{n+r}(F) \text{ for all } u\in U  \}
\end{align*}
and call it the $\Omega$-annihilator resp. the $\nu$-annihilator of $U$.
\end{defi}
Note that $\ann \nu^n_F(U) = \ann \Omega^n_F(U) \cap \nu_n(F)$ for every $U \subset \Omega^r(F)$. For non-empty sets $S_1,\ldots, S_r \subset F$ we set 
$\de S_1 \w \ldots \w \de S_r := \{ \de s_1 \w \ldots \w \de s_r \mid s_1 \in S_1, \ldots, s_r \in S_r \}$ and replace the slot $\de S_i$ by $\de s_i$ if 
$S_i=\{s_i\}$. The following lemma is now easy to check.

\begin{lem}{\cite[Corollary 2.3]{q58}}\label{33}
Let $S_1,\ldots, S_r \subset F$ be non-empty sets with $\pdeg_F (S_i)=k_i \in \NN$ for $i = 1,\ldots r$ and 
choose a $p$-basis $\{a_{i1},\ldots, a_{ik_i}\}$ of $F^p(S_i)$. Then
\begin{align*}
\ann\Omega^n_F \left( \de S_1 \w \ldots \w \de S_r \right) &= \ann\Omega^n_F \left( \de \,\{a_{11},\ldots,a_{1k_1}\} \w \ldots \w \de\, \{a_{r1},\ldots,a_{rk_r}\} \right)\\
&=\ann\Omega^n_F \left( \de \left( F^p(S_1)\right) \w \ldots \w \de \left( F^p(S_r)\right) \right).
\end{align*}
\end{lem}

Specific annihilators are rather difficult to compute. The known general cases are listed in the following three propositions.

\begin{pro}{\cite[Proposition 2.4]{q58}}\label{34}
For non-empty sets $S_1,\ldots,S_r \subset F$ with $\pdeg_F(S_i)=k_i\in \NN$ and $p$-bases $\{a_{i1},\ldots,a_{ik_i}\}$ of $F^p(S_i)$ 
for $i=1,\ldots,r$ say 
$\pdeg (S_1 \cup \ldots \cup S_r)=k_1+ \ldots +k_r$. Then 
\begin{align*}
\ann\Omega^n_F \left( \de S_1 \w \ldots \w \de S_r \right)= 
\sum_{i=1}^r \de a_{i1}\w\ldots\w\de a_{ik_i} \w \Omega^{n-k_i}(F).
\end{align*}
\end{pro}

\begin{pro}{\cite[Proposition 2.5]{q58}}\label{35}
Let $S \subset F$ be a non-empty set with $\pdeg_F(S)=k\in\NN$ and let $\{a_1,\ldots, a_k\}$ be a $p$-basis of $F^p(S)$. For $r\in \{1,\ldots,k\}$ set 
$t=k-r+1$. Then 
\begin{enumerate}[(a)]
	\item \begin{align*}
	\ann\Omega^n_F\left(\bigwedge\nolimits^r \de S \right) =\sum_{\substack{\{i_1,\ldots,i_t\}\subset\{1,\ldots,k\} \\ i_1<\ldots <i_t  }} 
	\de a_{i_1}\w\ldots \w\de a_{i_t} \w\Omega^{n-t}(F)
\end{align*}
and for $r>k$ we have $	\ann\Omega^n_F\left(\bigwedge\nolimits^r \de S \right)=\Omega^n(F)$.

	\item if we additionally assume $F^{p-1}=F$, we have
\begin{align*}
	\ann\nu^n_F\left({\bigwedge}^r\de S \right) = \bigg[ \frac{\de y_1}{y_1}\w\ldots\w\frac{\de y_t}{y_t} ~\Big\vert~ y_1,\ldots,y_t\in F^p(a_1,\ldots,a_k)^*  \bigg] \w\nu_{n-t}(F)
		\end{align*}
	and for $r>k$ we have $\ann\nu^n_F({\bigwedge}^r\de S)=\nu_n(F)$.
\end{enumerate}
\end{pro}

\begin{pro}{\cite[Proposition 2.6]{q58}}\label{36}
Let $S_1,\ldots,S_{r+1}\subset F$ be non-empty sets with $\pdeg_F (S_i)=1$ for $i=1,\ldots,r$ and $\de S_1 \w \ldots \de S_r \w \de S_{r+1}\neq \{0\}$. 
Set $\pdeg_F(S_1 \cup \ldots \cup S_r \cup S_{r+1})=r+\ell$ with $\ell \in \ZZ$. Then $\ell \geq 1$. Now choose a $p$-independent set 
$\{a_1,\ldots,a_r,e_1,\ldots, e_{\ell}\}\subset F$ such that $F^p(S_i)=F^p(a_i)$ for $i=1,\ldots,r$ 
and $F^p(a_1,\ldots, a_r)(S_{r+1})=F^p(a_1,\ldots,a_r)(e_1,\ldots,e_{\ell})$. Then
\begin{enumerate}[(a)]
	\item we have \begin{align*}
	\ann\Omega^n_F( \de S_1 \w \ldots \w \de S_r\w \de S_{r+1} ) = \sum_{i=1}^r \de a_i \w \Omega^{n-1}(F) + \de e_1\w\ldots \w\de e_{\ell}\w\Omega^{n-\ell}(F).
\end{align*}
	
	\item if we additionally assume $F^{p-1}=F$, we have
		\begin{align*}
		&\ann\nu^n_F(\de S_1 \w \ldots \w \de S_r \w \de S_{r+1}) =\bigg[ \frac{\de x}{x}  ~\Big\vert~ x\in F^p(a_1,\ldots,a_r)^* \bigg]\w \nu_{n-1}(F) \\
		&\ \ + \bigg[\frac{\de y_1}{y_1}\w\ldots\w\frac{\de y_{\ell}}{y_{\ell}} ~\Big\vert~ 
		y_1,\ldots,y_{\ell}\in F^p(a_1,\ldots,a_r,e_1,\ldots e_{\ell})^* \bigg]\w\nu_{n-\ell}(F).
		\end{align*}
\end{enumerate}
\end{pro}

\begin{rem}\label{37}
The assumption $F^{p-1}=F$ in the Propositions \ref{35} and \ref{36} seem to come out of nowhere, but are needed here since 
the proofs rely on the use of the Lemma of Kato \cite[Lemma 2]{q4} which also assumes this property. It is open (but likely) that 
these characterizations of the annihilators are still correct if the assumption  $F^{p-1}=F$ is dropped. In later sections, we will 
mainly use these results for the case $p=2$ where this assumption becomes redundant.
\end{rem}

We will end this section with an interesting connection between $\Omega$-kernels of modular purely inseparable extensions and annihilators, which we will need in later sections.

\begin{pro}{\cite[Proposition 4.3]{q58}}\label{38}
Let $\{ b_1,\ldots, b_r \}\subset F$ be $p$-independent, $m_1,\ldots,m_r\in\NN$ and set $E=F\left(\sqrt[p^{m_1}]{b_1}, \ldots, \sqrt[p^{m_r}]{b_r} \right)$. 
Then 
\begin{align*}
\Omega^n(E/F)&= \sum_{i=1}^r \de b_i \w \Omega^{n-1}(F)= \ann\Omega^n_F(\de b_1 \w \ldots \w \de b_r).
\end{align*}
\end{pro}

\end{subsection}

\end{section}

\begin{section}{Differential forms under iterated function field extensions of $p$-forms}\label{ffe}

In this section, we will combine the previous results to rewrite the $\Omega$-kernels of iterated function field extensions of $p$-forms 
as an annihilator of specific sets related to the $p$-forms. For $p$-forms $\p_1,\ldots,p_r$ defined over $F$, the iterated function field  
$F(\p_1,\ldots,\p_r)$
is defined as the field $F_r$ using the recursion $F_0:=F$ and $F_i:=F_{i-1}(\p_i)$. To simplify notations and computations in the upcoming 
results, let us fix some assumptions on the forms $\p_1,\ldots,\p_r$ which will not pose any restrictions on the next results.

\begin{rem}\label{41}
For $i=1,\ldots, r$ let $\p_i$ be a $p$-form over $F$ of dimension $1+n_i$ with $n_i\in\NN$. After scaling, 
we may assume $\p_i= \lr 1, a_{i1}, \ldots, a_{in_i}\rr _p$ with $a_{i1}, \ldots, a_{in_i}\in F$ (note that by Lemma \ref{23} this scaling does not change the 
field extension). Set $\ndeg_F(\p_i)=p^{k_i}$ and assume that the $a_{i1}, \ldots, a_{in_i}$ are ordered such that 
$N_F(\p_i)=F^p(a_{i1}, \ldots, a_{ik_i})$. Let $X_{i1},\ldots,X_{in_i}$ be variables for $i=1,\ldots,r$ and set $\p_i = \lr 1 \rr_p \perp \p_i'$, i.e. 
$\p_i'=\lr  a_{i1}, \ldots, a_{in_i}\rr _p$. We define $T_i=\p_i'(X_{i1},\ldots,X_{in_i})\in F[X_{ij} \mid j=1,\ldots,n_i]$. With these notations we have 
$$F(\p_1,\ldots,\p_r) \cong F\left(X_{ij} \mid i=1,\ldots,r\ , \ j=1, \ldots, n_i \right)\left( \sqrt[p]{T_1},\ldots, \sqrt[p]{T_r} \right).$$
Now set $M:=F\left(X_{ij} \mid i=1,\ldots,r\ , \ j=1, \ldots, n_i \right)$ and note that this representation of $F(\p_1,\ldots,\p_r)$ shows that 
it is independent of the ordering of the $\p_1,\ldots,\p_r$. Since $M\left( \sqrt[p]{T_1},\ldots, \sqrt[p]{T_r} \right)/M$ is a purely inseparable extension of 
exponent 1, after possibly reordering the $T_1,\ldots, T_r$ (and with it the $\p_1,\ldots,\p_r$), we find $s\leq r$ with 
$[M\left( \sqrt[p]{T_1},\ldots, \sqrt[p]{T_r} \right):M]=p^s$ and $M\left( \sqrt[p]{T_1},\ldots, \sqrt[p]{T_r} \right) = 
M\left( \sqrt[p]{T_1},\ldots, \sqrt[p]{T_s} \right)$. Set $L:=F\left(X_{ij} \mid i=1,\ldots,s\ , \ j=1, \ldots, n_i \right)$. Now it is easy to see that 
the extension $F(\p_1,\ldots,\p_r)/F(\p_1,\ldots,\p_s)$ is purely transcendental and thus we have 
$$\Omega^n(F(\p_1,\ldots,\p_r)/F) =\Omega^n(F(\p_1,\ldots,\p_s)/F)$$
by Lemma \ref{315}.
\end{rem}

Before we are able to proof our main theorem of this section, we will start  with some useful lemmas.

\begin{lem}\label{42}
Let $\p_1,\ldots,\p_r$ be $p$-forms over $F$ as in Remark \ref{41}. For $i\in\{1,\ldots,r\}$ we find $q_{ij}\in F[X_{ij} \mid j=1,\ldots, n_i]$ with 
$$\de T_i = q_{i1} \de a_{i1} + \ldots +q_{ik_i} \de a_{ik_i}.$$
For a fixed $j \in \{1,\ldots,k_i\}$ there further exists a substitution such that 
$$q_{ij} \mapsto 1 \text{ and } q_{im} \mapsto 0 \text{ for all } m\in\{1,\ldots,n_i\}\setminus \{j\}.$$
In particular there exists a substitution $\de T_i \mapsto \de a_{ij}$ for all $j\in\{1,\ldots,k_i\}.$
\end{lem}

\begin{prf}
Fix $i\in\{1,\ldots,r\}$. By $T_i = \p_i'(X_{i1},\ldots,X_{in_i})=a_{i1}X_{i1}^p + \ldots + a_{in_i}X_{in_i}^p$ we have 
\begin{equation}\label{e1}
\de T_i = X_{i1}^p \de a_{i1} + \ldots + X_{in_i}^p \de a_{in_i}.
\end{equation}
Since $\{a_{i1},\ldots,a_{ik_i}\}$ is a $p$-basis of $F^p(a_{i1},\ldots,a_{in_i})$, by Lemma \ref{31} for $t=k_i+1,\ldots,n_i$ we find a representation 
$$\de a_{it} = \sum_{j=1}^{k_i} \lambda_{itj} \de a_{ij} \ \text{ with }\ \lambda_{itj}\in F.$$
Inserting this into Equation \eqref{e1} gives us 
\begin{align*}
	\de T_{i} = \sum_{j=1}^{k_i} X_{ij}^p \de a_{ij} + \sum_{t=k_i+1}^{n_i} \left( \sum_{j=1}^{k_i} \lambda_{itj} \de a_{ij} \right)   X_{it}^p  
		= \sum_{j=1}^{k_i} \left(  X_{ij}^p  + \sum_{t=k_i+1}^{n_i}  \lambda_{itj}  X_{it}^p \right) \de a_{ij}.
	\end{align*}
So setting $q_{ij}:= X_{ij}^p  + \sum_{t=k_i+1}^{n_i}  \lambda_{itj}  X_{it}^p$ gives us the representation we were looking for and 
for fixed $j\in\{1,\ldots,k_i\}$, the substitution $X_{ij} \mapsto 1$ and $X_{im}\mapsto 0$ for $m\in \{1,\ldots,n_i\}\setminus \{j\}$ produces
$q_{ij}\mapsto 1$ and $q_{ih}\mapsto 0$ for $m\in \{1,\ldots,k_i\}\setminus \{j\}$.\be 
\end{prf}

It is a natural question to ask what conditions for the $\p_1,\ldots, \p_r$ 
are necessary so that $\{T_1,\ldots,T_r\}$ is $p$-independent over $M$ (i.e. $s=r$ with the notations form Remark \ref{41}). 
An answer to this question gives the following lemma.

\begin{lem}\label{43}
Let $\p_1,\ldots,\p_r$ be $p$-forms over $F$ as in Remark \ref{41}. The following are equivalent
\begin{enumerate}[(a)]
\item $\{T_1,\ldots,T_r\}$ is $p$-independent over $M$.
\item $\de\, \{a_{11},\ldots, a_{1k_1}\} \w \ldots \w \de\, \{a_{r1},\ldots,a_{rk_r}\} \neq \{0\}.$
\item $\de N_F(\p_1) \w \ldots \w \de N_F(\p_r) \neq \{0\}.$
\end{enumerate}
\end{lem}

\begin{prf}

Let us start with the equivalence of (a) and (b). By Lemma \ref{42} we have 
	\begin{align}\label{e2}
	 \de T_1 \w \ldots \w \de T_r \overset{\ref{42}}{=}\sum_{j_1=1}^{k_1}\ldots\sum_{j_r=1}^{k_r} q_{1j_1}\ldots q_{rj_r} \de a_{1j_1}\w\ldots\w \de a_{rj_r}.
	\end{align}
Let us now assume (a) to be false. So let $\{T_1,\ldots,T_r\}$ be $p$-dependent over $M$ and thus $\de T_1 \w \ldots \w \de T_r =0$ by Lemma \ref{31}. Now using Lemma \ref{42} once again, for every $(j_1,\ldots,j_r)\in \bigtimes_{i=1}^r \{1,\ldots,k_i\}$ we can find a suitable substitution in 
Equation \eqref{e2} such that  
$ q_{1j_1}\mapsto 1,  \ldots , q_{rj_r}\mapsto 1$ and $q_{ij} \mapsto 0$ for all other $q_{ij}$. This gives us 
$$0=\de a_{1j_1}\w\ldots\w \de a_{rj_r} \ \text{ for all }\ (j_1,\ldots,j_r)\in \bigtimes_{i=1}^r \{1,\ldots,k_i\}$$
which then implies the negation of (b).

Now assume (b) to be false, i.e. $\de\, \{a_{11},\ldots, a_{1k_1}\} \w \ldots \w \de\, \{a_{r1},\ldots,a_{rk_r}\} =\{0\}.$ Thus the right hand side 
of Equation \eqref{e2} is $0$. Hence $\de T_1 \w \ldots \w \de T_r=0$ over $M$ which gives us the $p$-dependence of $\{T_1,\ldots,T_r\}$ over $M$ by Lemma \ref{31}. So now we have proven 
(a)$\iff$(b). The equivalence of (b) and (c) now easily follows by Lemma \ref{33}. \be
\end{prf}

\begin{lem}\label{44}
Let $\p_1,\ldots,\p_r$ be $p$-forms over $F$ as in Remark \ref{41} and let $\omega\in \Omega^n(F)$. The following are equivalent
\begin{enumerate}[(a)]
\item $\omega \in \Omega^n(F(\p_1,\ldots,\p_r) /F)$.
\item $\omega_M \in \Omega^n(M\left( \sqrt[p]{T_1},\ldots,\sqrt[p]{T_r}\right) /M)$.
\item $\omega_L \in \Omega^n(L\left( \sqrt[p]{T_1},\ldots,\sqrt[p]{T_s}\right) /L)$.
\end{enumerate}
\end{lem}

\begin{prf}
This can be easily checked by the definitions of the fields and is left to the reader. \be 
\end{prf}

Now we have assembled all necessary tools to prove the main theorem of this section.

\begin{theo}\label{45}
Let $\p_1,\ldots,\p_r$ be $p$-forms over $F$ as in Remark \ref{41}. Then 
$$\de \left\{ a_{11},\ldots,a_{1k_1}\right\} \w \ldots \w \de \left\{a_{s1},\ldots,a_{sk_s} \right\} \neq \{0\}$$
and we have 
	\begin{align*}
	\Omega^n(F(\p_1,\ldots,\p_r)/F) &=\Omega^n( F(\p_1,\ldots,\p_s)/F) \\
	&= \ann\Omega^n_F (\de \left\{ a_{11},\ldots,a_{1k_1}\right\} \w \ldots \w \de \left\{a_{s1},\ldots,a_{sk_s} \right\}) \\
	&=\ann\Omega_F^n( \de N_F(\p_1) \w \ldots \w \de  N_F(\p_s))
	\end{align*}
as well as
	\begin{align*}
	\nu_n(F(\p_1,\ldots,\p_r)/F) &= \nu_n(F(\p_1,\ldots,\p_s)/F) \\
	&= \ann\nu^n_F (\de \left\{ a_{11},\ldots,a_{1k_1}\right\} \w \ldots \w \de \left\{a_{s1},\ldots,a_{sk_s} \right\}) \\
	&=\ann\nu_F^n( \de N_F(\p_1) \w \ldots \w \de  N_F(\p_s)).
	\end{align*}
In particular the kernels $\Omega^n(F(\p_1,\ldots,\p_r)/F)$ and $\nu^n(F(\p_1,\ldots,\p_r)/F)$ only depend on the norm fields of the forms $\p_1,\ldots,\p_s$ but 
not on the forms itself.
\end{theo}

\begin{prf}
We have $\de \left\{ a_{11},\ldots,a_{1k_1}\right\} \w \ldots \w \de \left\{a_{s1},\ldots,a_{sk_s} \right\} \neq \{0\}$ by Lemma \ref{43} and the 
$p$-independence of $\{T_1,\ldots,T_s\}$ over $M$ (and over $L$) by construction. It is now sufficient to proof the first chain of equalities since the second one can 
be produced by the intersection with $\nu_n(F)$.

Since $F(\p_1,\ldots,\p_r)/ F(\p_1,\ldots,\p_s)$ is purely trancendental, any form $\omega\in \Omega^n(F)$ which becomes the zero form over 
$F(\p_1,\ldots,\p_r)$ must already be the zero form over $F(\p_1,\ldots,\p_s)$ by Lemma \ref{315}. Thus we have the equality $\Omega^n(F(\p_1,\ldots,\p_r)/F) 
=\Omega^n( F(\p_1,\ldots,\p_s)/F)$. 

The equality of the annihilators $\ann\Omega^n_F (\de \left\{ a_{11},\ldots,a_{1k_1}\right\} \w \ldots \w \de \left\{a_{s1},\ldots,a_{sk_s} \right\})$ 
and $\ann\Omega_F^n( \de N_F(\p_1) \w \ldots \w \de  N_F(\p_s))$ again directly follows by Lemma \ref{33}.

So all that is left to prove is the middle equation. Let us start with the inclusion $(\subseteq)$. Let 
$\omega\in \Omega^n(F(\p_1,\ldots,\p_s)/F)$. By Lemma \ref{44} we have $\omega_L \in \Omega^n(L\left( \sqrt[p]{T_1},\ldots,\sqrt[p]{T_s}\right) /L)$ and by
Proposition \ref{38} we may rewrite the $\Omega$-kernel $\Omega^n(L\left( \sqrt[p]{T_1},\ldots,\sqrt[p]{T_s}\right) /L)$ as the annihilator $\ann\Omega^n_L(\de T_1\w \ldots \w \de T_s)$ 
since the set $\{T_1,\ldots,T_s\}$ is $p$-independent over $L$. This gives us 
$\omega_L \w \de T_1\w \ldots \w \de T_s =0$ over $L$. Now again by Lemma \ref{42} we find $q_{ij}\in F[X_{ij} \mid j=1,\ldots,n_i]$ for $i=1,\ldots,s$ with 
	\begin{align}\label{e3}
	\omega_L \w \de T_1 \w \ldots \w \de T_s = \sum_{j_1=1}^{k_1}\ldots\sum_{j_s=1}^{k_s}q_{1j_1}\ldots q_{sj_s}\omega_L \w \de a_{1j_1}\w\ldots\w \de a_{sj_s}.
	\end{align}
Now using a suitable substitution in \eqref{e3} and using the fact that the left hand side of \eqref{e3} is zero, we get 
	\begin{align*}
	0= \omega \w \de a_{1j_1}\w\ldots\w \de a_{sj_s} \ \text{ for all } (j_1,\ldots,j_s)\in\bigtimes_{i=1}^s \{1,\ldots,k_i\},
	\end{align*}
and hence $\omega\in \ann\Omega^n_F( \de \left\{ a_{11},\ldots,a_{1k_1}\right\} \w \ldots \w \de \left\{a_{s1},\ldots,a_{sk_s} \right\})$.

For the reverse inclusion $(\supseteq)$ choose  $\omega \in \Omega^n(F)$ with 
 $\omega \w \de a_{1j_1}\w\ldots\w \de a_{sj_s}=0$ for all $(j_1,\ldots,j_s)\in\bigtimes_{i=1}^s \{1,\ldots,k_i\}$. Inserting this in 
 Equation \eqref{e3}, we readily get $\omega_L \w \de T_1 \w \ldots \w \de T_s=0$ and thus $ \omega_L \in \ann\Omega^n_L(\de T_1\w \ldots \w \de T_s)$. 
 This proof is then finished by rewriting this annihilator as a kernel using 
 Proposition \ref{38} and applying Lemma \ref{44} one last time. \be

\end{prf}

Since we now have established a connection between $\Omega$-kernels and $\nu$-kernels of iterated function field extensions of $p$-forms and 
annihilators. This enables us to apply our results \ref{34}, \ref{35} and \ref{36} to describe some new kernels.

\begin{theo}\label{46}
Let $\p_1,\ldots,\p_r$ be $p$-forms over $F$ as in Remark \ref{41}. 
\begin{enumerate}[(a)]
	\item If $\pdeg_F ( \{ a_{ij} \mid i=1,\ldots,s$, $j=1,\ldots,k_i \})= k_1+\ldots +k_s$, then we have
		$$ \Omega^n(F(\p_1,\ldots,\p_r)/F)= \sum_{i=1}^s \Omega^n(F(\p_i)/F)
		=\sum_{i=1}^s \de a_{i1}\w \ldots \w \de a_{ik_i} \w \Omega^{n-k_i}(F).$$

	\item Assume $\p_1=\ldots =\p_r=\p$ with $\ndeg_F (\p)=p^k$ and $N_F(\p)=F^p(a_1,\ldots,a_k)$ for $a_1,\ldots,a_k\in F$. Now set $t:=k-r+1$, then we have
		\begin{align*}
			\Omega^n(F(\underbrace{\p,\ldots ,\p}_{r \text{ times}})/F)
				= \begin{cases} \sum\limits_{i=1}^k \de a_i \w \Omega^{n-1}(F)  &,r\geq k, \\
				\sum\limits_{\substack{\{i_1,\ldots,i_t\} \subset \{1,\ldots,k\} \\ i_1<\ldots < i_t}} \de a_{i_1}\w\ldots\w\de a_{i_t}\w\Omega^{n-t}(F) \ \   &,r<k	.			
				\end{cases}
		\end{align*}
		If we additionally assume $F=F^{p-1}$, then 
		\begin{align*}
			\nu_n(F(\underbrace{\p,\ldots ,\p}_{r \text{ times}})/F)
				= \begin{cases}\left[ \frac{\de x}{x} ~\big\vert~ x\in N_F(\p)^* \right]\w\nu_{n-1}(F)  &,r\geq k, \\
				\left[ \frac{\de y_1}{y_1}\w\ldots\w \frac{\de y_t}{y_t} ~\big\vert~ y_1,\ldots,y_t \in N_F(\p)^* \right] \w \nu_{n-t}(F) \ \   &,r<k.				
				\end{cases}
		\end{align*}
	\item Say $\ndeg_F (\p_1)=\ldots = \ndeg_F (\p_{s-1})=p$ and set $a_{11}=a_1,\ldots, a_{s-1,1}=a_{s-1}$. Then $\{a_1,\ldots,a_{s-1}\}$ is $p$-independent over $F$. 
	Let $\pdeg_F ( \{ a_{1},\ldots,a_{s-1}\} \cup \{ a_{s1},\ldots,a_{sk_s}\})= s-1+\ell$ and choose suitable $e_1,\ldots,e_{\ell}\in F$ such that  
	$F^p(a_{1},\ldots,a_{s-1})(a_{s1},\ldots,a_{sk_s})=F^p(a_{1},\ldots,a_{s-1})(e_1,\ldots,e_{\ell})$.
	Then $\ell \geq 1$ and we have 
			$$\Omega^n(F(\p_1,\ldots,\p_r)/F)= \sum_{i=1}^{s-1} \de a_{i} \w \Omega^{n-1}(F) + \de e_1 \w \ldots \w \de e_{\ell} \w\Omega^{n-\ell}(F).$$
		If we additionally assume $F=F^{p-1}$, then
			\begin{align*}
				&\nu_n (F(\p_1,\ldots,\p_r)/F) = \left[ \frac{\de x}{x} ~\Big\vert~ x\in F^p(a_{1},\ldots,a_{s-1})^* \right] \w \nu_{n-1}(F) \\
				&\ \ + \left[ \frac{\de y_1}{y_1}\w \ldots\w\frac{\de y_{\ell}}{y_{\ell}} ~\Big\vert~ y_1,\ldots,y_{\ell}\in F^p(a_{1},\ldots a_{s-1},e_1,\ldots,e_{\ell})
				 \right] \w \nu_{n-\ell}(F).
			\end{align*}
			\end{enumerate} 
\end{theo}


We can now use Theorem \ref{46}(b) to describe the $\Omega$-kernel and the 
$\nu$-kernel of the compositum of a multiquadratic purely inseparable extension of exponent 1 with a function field extension of a $p$-form. To do so, recall that for a $p$-form $\p$ over $F$ with $\ndeg_F(\p)=p$, we can find $a\in F$ with $N_F(\p)=F^p(a)$. Then $F(\p)$ is a purely transcendental 
extension of $F(\sqrt[p]{a})$ and thus $\Omega^n(F(\p)/F)=\Omega^n(F(\sqrt[p]{a})/F)$ (To be more precise:  $F(\sqrt[p]{a})$ is the algebraic closure of $F$ in $F(\p)$, 
see \cite[Proposition 7.6]{q11}). Hence, in the situation of Theorem \ref{46}(c) the extension 
$$F(\p_1,\ldots,\p_s) / F\left(\sqrt[p]{a_1}, \ldots, \sqrt[p]{a_{s-1}}\right)(\p_s) $$
is purely transcendental. Note that $\ndeg_{F(\p_1,\ldots,\p_{i-1})}(\p_i)=p$ and $(F(\p_1,\ldots,\p_{i-1}))^p(a_{i})=N_{F(\p_1,\ldots,\p_{i-1})}(\p_i)$ 
still hold for all $i\in\{1,,\ldots,s-1\}$ since $\de T_1 \w \ldots \w \de T_{s-1}\neq 0$. So now we get

\begin{theo}\label{47}
Let $\{a_1,\ldots,a_s\}\subset F$ be $p$-independent and let $\p=\lr 1,b_1,\ldots,b_t\rr_p$ be a (not necessarily anisotropic) $p$-form over $F$ with 
$\ndeg_{F(\sqrt[p]{a_1},\ldots,\sqrt[p]{a_s})}(\p)>1$. Say $\pdeg_F(\{a_1,\ldots,a_s\}\cup \{b_1,\ldots,b_t\})=s+\ell$ with 
$F^p(a_1,\ldots,a_s)(b_1,\ldots,b_t)=F^p(a_1,\ldots,a_s)(e_1,\ldots,e_{\ell})$ for some suitable $e_1,\ldots,e_{\ell}\in F$. Then 
$\ell \geq 1$ and we have
\begin{align*}
	\Omega^n(F(\sqrt[\leftroot{2}\uproot{2} p]{a_1},\ldots,\sqrt[\leftroot{2}\uproot{2} p]{a_s})(\p)/F)&= \ann \Omega^n_F (\de a_1 \w \ldots \w \de a_s \w \de N_F(\p)) \\
	&=\sum_{i=1}^s \de a_i \w \Omega^n(F) + \de e_1\w \ldots \w \de e_{\ell} \w \Omega^{n-\ell}(F)
	\end{align*}
and 
	$$ \nu_n(F(\sqrt[\leftroot{2}\uproot{2} p]{a_1},\ldots,\sqrt[\leftroot{2}\uproot{2} p]{a_s})(\p)/F)= \ann \nu^n_F (\de a_1 \w \ldots \w \de a_s \w \de N_F(\p)). $$
If we additionally assume $F=F^{p-1}$, then  
	\begin{align*}
		\nu_n(&F(\sqrt[\leftroot{2}\uproot{2} p]{a_1},\ldots,\sqrt[\leftroot{2}\uproot{2} p]{a_s})(\p)/F)= \left[\frac{\de x}{x} ~\Big\vert~ x\in F^p(a_1,\ldots,a_s)^* \right]\w \nu_{n-1}(F) \\
			&+\left[ \frac{\de y_1}{y_1}\w\ldots\w\frac{\de y_{\ell}}{y_{\ell}} ~\Big\vert~ 
			y_1,\ldots,y_{\ell}\in F^p(a_1,\ldots,a_s,e_1,\ldots e_{\ell}) \right] \w \nu_{n-\ell}(F).
	\end{align*}

\end{theo}

\end{section}

\begin{section}{Transfer to the theory of bilinear forms in characteristic two}\label{bf}

In this last section we will translate the previous results into the language of bilinear forms over fields of characteristic two. So for that, from now on let 
$F$ be a field of characteristic two. We will briefly recap the main aspects of the theory, but we highly suggest to look up further details, for example 
in \cite{q14} or in \cite{q16}.

By a bilinear form over $F$, we always mean a finite dimensional symmetric non-degenerate bilinear form over $F$. Now let $\bb: V \times V \to F$ be a 
bilinear form on the vector space $V$. We call $\bb$ isotropic if there is a $v \in V \setminus \{0\}$ with $\bb(v,v)=0$ and anisotropic otherwise. A subspace 
$W\subset V$ is called isotropic, if $\bb(w,w)=0$ for all $w \in W$ and the form $\bb$ is called metabolic, if there exists an isotropic subspace $W\subset V$ 
with $2\dim W= \dim V$. Now $D_F(\bb):=\{\bb(v,v)\neq 0 \mid v \in V\}$ is called the set of elements represented by $\bb$. The bilinear form $\bb$ can be diagonalized iff 
$D_F(\bb)\neq \emptyset$ and in this case we write $\bb=\lr a_1,\ldots,a_n\rr$ with $a_1,\ldots,a_n \in F$.

Recall that every bilinear form has a orthogonal decomposition into an anisotropic and a metabolic part. The anisotropic part is determined uniquely up to 
isometry and is called $\bb_{an}$. Two bilinear forms $\bb$ and $\bb'$ are called Witt-equivalent if $\bb_{an}\cong \bb'_{an}$. The classes of 
Witt equivalent forms then form a ring $W(F)$, called the Witt ring of $F$. The addition in $W(F)$ is defined using the orthogonal sum of forms and the 
multiplication is induced by the tensor product of forms. 
By abuse of notation we will write $\bb$ for both the bilinear form and 
Witt class of $\bb$. 

A given field extension $L/F$ induces natural homomorphisms $\iota : W(F) \to W(L)$ whose 
kernel $W(L/F)$ is called the bilinear Witt kernel
of the extension $L/F$. For $\bb \in W(F)$ viewed as an element of $W(L)$ we shortly write $\bb_L$. Note that $W(L/F)=\{ \bb\in W(F) \mid \bb_L \text{ is metabolic} \}$.

A bilinear form of the type $\lr 1,a_1\rr \otimes \ldots \otimes \lr 1,a_n\rr=:\lrp a_1,\ldots,a_n\rrr$ is called an $n$-fold bilinear Pfister form. 
A well known fact states that bilinear Pfister forms are either anisotropic or metabolic. Note that $\lrp a_1,\ldots,a_n\rrr$ is anisotropic if and only if 
the set $\{a_1,\ldots,a_n\}$ is $2$-independent over $F$. The Witt classes of even dimensional bilinear forms form an ideal $I(F)$ called the fundamental ideal 
of $F$. Its powers $I^n(F):=(I(F))^n$ are generated (as an additive group as well as an ideal) by the set of bilinear $n$-fold Pfister forms for all $n\in\NN$. 
Due to the inclusion chain 
$$W(F) \supset I(F) \supset I^2(F)\supset \ldots$$
the quotients $\overline{I^n(F)}:=I^n(F) / I^{n+1}(F)$ give rise to the graded Witt ring.

For a non-empty set $U\subset W(F)$ we call 
$$\ann\bb_F(U):= \{ \bb \in W(F) \mid \bb \otimes u =0 \text{ for all } u \in U\}$$
the bilinear annihilator of $U$. For non-empty sets $S_1,\ldots,S_r \subset F^*$ we further define 
$$ \lrp S_1\rrr \otimes \ldots \otimes \lrp S_r \rrr:= \lrp S_1,\ldots,S_r\rrr:= \{\lrp s_1,\ldots,s_r\rrr \mid s_1 \in S_1, \ldots, s_r \in S_r\}.$$

Note that the map $\pb: V \to F, v \mapsto \pb(v):=\bb(v,v)$ is a 2-form in the sense of section \ref{pform}, thus all results and definitions from section 
\ref{pform} apply to $\pb$. The form $\pb$ is called the polar form of $\bb$. The function field of $\bb$ is then defined as $F(\pb)$.

With these notations we are able to introduce the first half of the Theorem of Kato, which gives an answer to Milnor 
conjecture in characteristic two and can also be used 
to transfer our results from the previous sections to the theory of bilinear forms.

\begin{theo}\label{51}{\cite[p. 494]{q9}}
For every $n\in\NN$ there is an isomorphism
\begin{align*}
e_n : I^n(F) / I^{n+1}(F) &\to \nu_n(F) \\
\lrp a_1,\ldots , a_n \rrr \!\!\!\! \mod I^{n+1}(F) &\mapsto   \frac{\de a_1}{a_1} \w \ldots \w \frac{\de a_n}{a_n}   
\end{align*}
\end{theo}

Now to translate the results from section \ref{ffe} to bilinear forms, we will use a standard procedure which goes back to the case of $F$ having a finite $p$-basis. 
A detailed description of the necessary steps can for example be found in \cite[p.4]{q59} or in \cite[Theorem 5.3]{q55}, so  
we only state the results in this paper.

For a given set of 
bilinear forms $\bb_1, \ldots, \bb_r$, we can now apply the steps and the reordering described in Remark \ref{41}, so that the polar forms 
$\p_{\bb_1},\ldots,\p_{\bb_r}$  meet the properties stated in 
Remark \ref{41}. With this we are able to translate our main results 
from section \ref{ffe}.

\begin{theo}\label{52}
Let $\bb_1, \ldots, \bb_r$ be bilinear forms over $F$ with polar forms $\p_{\bb_1},\ldots,\p_{\bb_r}$ as described in Remark \ref{41}. Then 
	\begin{align*}
	W(F(\bb_1, \ldots, \bb_r)/F) &=W( F(\p_{\bb_1},\ldots,\p_{\bb_s})/F) \\
	&= \ann\bb_F (\lrp a_{1j_1},\ldots, a_{sj_s}\rrr \mid j_i=1,\ldots,k_i \text{ for } i=1,\ldots,s)\\
	&=\ann\bb_F( \lrp N_F(\p_{\bb_1}),  \ldots, N_F(\p_{\bb_s})\rrr).
	\end{align*}
\end{theo}

\begin{theo}\label{53}
Let $\bb_1, \ldots, \bb_r$ be bilinear forms over $F$ with polar forms $\p_{\bb_1},\ldots,\p_{\bb_r}$ as described in Remark \ref{41}.
\begin{enumerate}[(a)]
	\item Assume $\bb_1=\ldots =\bb_r=\bb$ with $\ndeg_F (\pb)=2^k$ and $N_F(\pb)=F^2(a_1,\ldots,a_k)$ for $a_1,\ldots,a_k\in F$. Now set $t:=k-r+1$, then we have
		\begin{align*}
			W(F(\underbrace{\pb,\ldots ,\pb}_{r \text{ times}})/F)
				= \begin{cases}\left[ \lrp x \rrr \mid x\in N_F(\pb)^* \right]  &,r\geq k, \\
				\left[ \lrp y_1,\ldots,y_t\rrr \mid y_1,\ldots,y_t \in N_F(\pb)^* \right]  \ \   &,r<k.			
				\end{cases}
		\end{align*}
		
	\item Say $\ndeg_F (\p_{\bb_1})=\ldots = \ndeg_F (\p_{\bb_{s-1}})=2$ and $\p_{\bb_s}=\lr 1,a_{s1}, \ldots,a_{sn_s}\rr_2$. 
	Set $N_F(\p_{\bb_i})=F^2(a_i)$ for $i=1,\ldots,s-1$. Then $\{a_1,\ldots,a_{s-1}\}$ is $2$-independent over $F$ and 
	say we have $2-\deg_F ( \{ a_{1},\ldots,a_{s-1}\} \cup \{ a_{s1},\ldots,a_{sk_s}\})= s-1+\ell$ with 
	$F^p(a_{1},\ldots,a_{s-1})(a_{s1},\ldots,a_{sk_s})=F^p(a_{1},\ldots,a_{s-1})(e_1,\ldots,e_{\ell})$ for some suitable $e_1,\ldots,e_{\ell}\in F$. 
	Then $\ell \geq 1$ and $W(F(\p_{\bb_1},\ldots,\p_{\bb_r})/F)$ is generated by the forms
	\begin{itemize}
	\item $ \lrp x\rrr$ with $x \in F^2( a_1,\ldots,a_{s-1})$,
	\item $ \lrp y_1,\ldots,y_{\ell}\rrr$ with $ y_1,\ldots,y_{\ell}\in F^2(a_{1},\ldots a_{s-1},e_1,\ldots,e_{\ell})$.
	\end{itemize}	
			\end{enumerate} 
\end{theo}


\begin{theo}\label{54}
Let $\{a_1,\ldots,a_s\}\subset F$ be $2$-independent, 
set $L= F\left(\sqrt{a_1},\ldots,\sqrt{a_s}\right)$ and let $\bb=\lr 1,b_1,\ldots,b_t\rr$ be a bilinear form over $F$ with $\ndeg_{L}(\pb)>1$. 
Say $2-\deg_F(\{a_1,\ldots,a_s\}\cup \{b_1,\ldots,b_t\})=s+\ell$ with 
$F^2(a_1,\ldots,a_s)(b_1,\ldots,b_t)=F^2(a_1,\ldots,a_s)(e_1,\ldots,e_{\ell})$ for some suitable $e_1,\ldots,e_{\ell}\in F$. Then 
$\ell \geq 1$ and the kernel $W\left(L(\pb)/F\right)$ is generated by the forms
	\begin{itemize}
	\item $ \lrp x\rrr$ with $x \in L^2$.
	\item $ \lrp y_1,\ldots,y_{\ell}\rrr$ with $ y_1,\ldots,y_{\ell}\in L^2(e_1,\ldots,e_{\ell})$.
	\end{itemize}
\end{theo}

\begin{rem}
\begin{enumerate}[(a)]
 
\item In \cite{q26} Aravire, Laghribi and O'Ryan computed the quadratic Witt kernel for a field extension of the form 
$F(\sqrt{a_1},\ldots,\sqrt{a_s})(\pi)/F$ with a bilinear Pfister form $\pi=\lrp c_1,\ldots,c_r\rrr$ such that $\{a_1,\ldots,a_s,c_1,\ldots,c_r\}$ is 
$2$-independent over $F$. So Theorem \ref{54} gives the bilinear counterpart to the results given in \cite{q26} for a more general field extension. Note that 
in Theorem \ref{54}, the set $\{a_1,\ldots,a_s,b_1,\ldots,b_t\}$ does not need to be $2$-independent and $\pi$ does not need to be a quasi-Pfister form. An interesting difference between these two kernels is 
the fact that in the quadratic case under the assumptions made in \cite{q26}, the kernel is the sum of the kernels for the extensions 
 $F(\sqrt{a_1},\ldots,\sqrt{a_s})/F$ and $F(\pi)/F$. In the bilinear case however, even with the stronger assumptions made in \cite{q26}, such an additive 
 behaviour for the kernels does not accrue.
 
 \item Note that in Theorem \ref{54}, the assumption of 
the $\{a_1,\ldots,a_s\}\subset F$ being $2$-independent over $F$ 
can always be assumed without restriction. This comes from the fact that 
every purely inseparable extension of exponent 1 is modular. To be more 
precise, if $E=F(\sqrt{c_1},\ldots, \sqrt{c_t})$ with $[E:F]=2^s$, we can 
always choose suitable $2$-independent $e_1,\ldots, e_s \in \{c_1,\ldots,c_t\}$ such that $E=F(\sqrt{e_1},\ldots, \sqrt{e_s})$.
 
 \item We are mainly interested in the bilinear Witt kernels of field extensions. But let us further define  
 a different kind of kernel. For a field extension $K/F$ and $m\in\NN$ set 
 $I^m(K/F)=\{ \bb \in I^m(F) \mid \bb_K=0\}=W(F)\cap I^m(F)$. This 
 set is obviously a group and the transfer from differential forms to 
 bilinear forms via Kato's Theorem also gives us additive generators 
 of the group $I^m(F(\bb_1, \ldots, \bb_r)/F)$ in the studied cases.
\end{enumerate}

\end{rem}

\end{section}

\bibliographystyle{abbrv}
\bibliography{literatur}

\end{document}